\documentclass[11pt,a4paper]{article}
\usepackage[T1]{fontenc}
\usepackage{lmodern}
\usepackage[margin=26mm]{geometry}
\usepackage{amsmath,amssymb,amsthm}
\usepackage{microtype}

\newtheorem{theorem}{Theorem}[section]
\newtheorem{lemma}[theorem]{Lemma}
\newtheorem{proposition}[theorem]{Proposition}
\theoremstyle{remark}

\newcommand{\F}{\mathcal F}
\newcommand{\ord}{\operatorname{ord}}
\newcommand{\supp}{\operatorname{supp}}
\newcommand{\s}{\mathsf s}
\newcommand{\T}{\mathsf T}
\newcommand{\h}{\mathsf h}
\newcommand{\mult}{\mathsf v}
\newcommand{\knum}{\mathsf k}
\newcommand{\wt}{\mathsf w}

\makeatletter
\newcommand{\blfootnote}[1]{%
  \begingroup
  \renewcommand{\thefootnote}{}%
  \renewcommand{\@makefntext}[1]{\noindent ##1}%
  \footnotetext{#1}%
  \addtocounter{footnote}{-1}%
  \endgroup
}
\makeatother

\numberwithin{equation}{section}

\title{Zero-sum subsequences of cross number $1$ over cyclic groups}
\author{
Weidong Gao\thanks{Corresponding author.} \\
Center for Applied Mathematics, Tianjin University \\
Tianjin 300072, P.R. China \\
E-mail: \texttt{wdgao@nankai.edu.cn}
\and
Yucen Mu \\
Center for Applied Mathematics, Tianjin University \\
Tianjin 300072, P.R. China \\
E-mail: \texttt{yucen4617@tju.edu.cn}
}
\date{}

\begin{document}
\maketitle
\begin{abstract}
Let $G$ be a finite abelian group, and let $\T(G)$ be the smallest integer $t$ such that every sequence of length $t$ over $G$ contains a nonempty zero-sum subsequence of cross number $1$. 
We determine the exact value of $\T(C_n)$ for every
positive integer $n$, showing that $\T(C_n)=2n-1$
when $n$ is odd and $\T(C_n)=3n/2-1$ when $n$ is even.
In particular, $\T(C_n)=\s(C_n)$ for odd $n$.
\end{abstract}


\section{Introduction and main results}

Given a positive integer $n$ and integers
$a_1,\ldots,a_{\ell}$, we ask whether there exists
a nonempty subset $I\subseteq[1,\ell]$ such that
$n\mid\sum_{i\in I}a_i$ and $\sum_{i\in I}\gcd(a_i,n)=n.$
We show that the sharp bound is $\ell=2n-1$ for odd $n$
and $\ell=3n/2-1$ for even $n$. This problem can be formulated
in terms of zero-sum subsequences of cross
number $1$ over cyclic groups.

Let $G$ be a finite abelian group, written additively. We denote by
$C_n$ a cyclic group of order $n$. For a general finite abelian group
$G$, we can decompose $G$ as a direct sum of cyclic groups
$C_{n_1}\oplus\cdots\oplus C_{n_r}$ such that
$1<n_1\mid\cdots\mid n_r\in\mathbb{N}$
(if $n_1=\cdots=n_r=n$, it will be abbreviated as $C_n^r$),
where $r$ and $n_r$ are respectively called the \emph{rank}
and \emph{exponent} of $G$.
Usually, the exponent of $G$ is simply denoted by $\exp(G)$.
The order of an element $g$ of $G$ will be written as
$\ord(g)$.

A sequence over $G$ is a finite unordered product
$S=g_1\cdot\ldots\cdot g_\ell$, where repetitions are allowed.
The cross number of $S=g_1\cdot\ldots\cdot g_\ell$ is
\[
\knum(S)=\sum_{i=1}^{\ell}\frac1{\ord(g_i)}.
\]

For a finite abelian group $G$ and a positive integer $k$,
let $\s_k(G)$ denote the smallest integer $t$ such that
every sequence over $G$ of length $t$ contains a
zero-sum subsequence of length $k\exp(G)$.
In particular, $\s_1(G)=\s(G)$ is the classical Erd\H{o}s--Ginzburg--Ziv
constant. The Erd\H{o}s--Ginzburg--Ziv theorem gives $\s(C_n)=2n-1$
for every positive integer $n$ (see \cite{[EGZ]}).
For further results on these constants,
we refer to~\cite{[GHPS],[GH2006]}.

Let $\T(G)$ be the smallest integer $t$ such that every sequence of length $t$ over $G$ contains a nonempty zero-sum subsequence of cross number $1$. In recent years, Gao, Hui, Jiang, Li and Wang introduced $\T(G)$ and studied its relation to classical zero-sum constants (see~\cite{[GHJLW]}). They determined $\T(C_n)$ for prime powers and established $\T(C_n)=\s(C_n)$ for several further classes of odd integers. They also conjectured that $\T(C_n^r)=\s(C_n^r)$ for all odd integers $n>1$ and positive integers $r$. We determine the invariant for all cyclic groups $C_n$.

\begin{theorem}\label{thm:main}
For every positive integer $n$,
\[
\T(C_n)=
\begin{cases}
2n-1,&n\text{ odd},\\[2pt]
\dfrac{3n}{2}-1,&n\text{ even}.
\end{cases}
\]
\end{theorem}

The lower bounds follow from the sequences in \cite[Lemma~3.1]{[GHJLW]}.
For $n=1$, there is nothing to prove. Let $n>1$ and let $g$ generate $C_n$.
If $n$ is odd, then $g^{n-1}\cdot(-g)^{n-1}$ contains no zero-sum
subsequence of length $n$. 
If $n$ is even, the same holds for
$g^{n-1}\cdot(-g)^{n/2-1}$. Since every term of these sequences has order
$n$, a subsequence $T$ satisfies $\knum(T)=|T|/n$. Hence, neither sequence
contains a zero-sum subsequence of cross number $1$.
For the upper bounds, we prove the stronger statement in Proposition~\ref{prop:sequence}.

The paper is organized as follows.
Some notation and some preliminary results are gathered in Section 2.
Sections~3 and~4 are devoted to the proofs of
Proposition~\ref{prop:sequence}(i) and (ii), respectively.

\section{Notation and preliminaries}

We denote by $\mathbb{N}$ the set of positive integers, by $\mathbb{N}_0$ the set $\mathbb{N} \cup \{0\}$ and by $\mathbb{Z}$ the set of integers.
For any $a,b \in \mathbb{Z}$, we set $[a,b]=\{x\in \mathbb{Z}\mid a\le x\le b\}$.
For any real number $x$, let $\lceil x \rceil$ be the smallest integer not less than $x$, and $\lfloor x \rfloor$ be the largest integer not greater than $x$.

Throughout this paper, let $p$ be a prime and $G$ be a finite additive abelian group with identity element $0$.
For a set $X$, we denote by $\mathcal{F}(X)$ the free (abelian, multiplicative)
monoid with basis $X$. The elements of $\mathcal{F}(X)$ are called \emph{sequences} over $X$.
We can write a sequence $S\in \mathcal{F}(X)$ in the form
$$
S=\prod_{g\in X}g^{\mult_g(S)},
$$
where $\mult_g(S)\in \mathbb{N}_0$ is called the \emph{multiplicity} of the element $g$ in $S$.
We say that $S$ \emph{contains} the element $g$ of $X$ if $\mult_g(S)>0$.
A sequence $S'$ is called a \emph{subsequence} of $S$, denoted by $S'\mid S$, if $\mult_g(S')\le \mult_g(S)$ for all $g\in X$.
The sequence $SS'^{-1}$ denotes the subsequence obtained from $S$ by deleting $S'$. Two subsequences $S'$ and $S''$ of $S$ are called \emph{disjoint} if $S'\mid SS''^{-1}$. The unit element $1 \in \mathcal{F}(X)$ is called the \emph{empty} sequence.

For a group $G$ and a sequence
$$S=g_1\cdots g_{\ell}=\prod_{g\in G}g^{\mult_g(S)}\in \mathcal{F}(G),$$
we define
\begin{itemize}
\item $|S| = \ell =\sum_{g\in G}\mult_{g}(S)\in \mathbb{N}_{0}$
the \emph{length} of $S$,

\item $\h(S) = \max\{\mult_g(S)\ |\ g \in G\}\in[0, |S|]$ the \emph{maximum of the multiplicities} of $S$,

\item $\sigma(S)=\sum_{i=1}^{\ell}g_{i}=\sum_{g\in G}\mult_{g}(S)g
\in G$ the \emph{sum} of $S$,

\item $\supp(S)=\{g\in G\mid \mult_g(S)>0\}\subseteq G$ the \emph{support} of $S$,

\item $S$ a \emph{zero-sum sequence} if $\sigma(S)=0\in G$.
\end{itemize}

Every map of abelian groups $\varphi:G\rightarrow H$ can be extended to a homomorphism $\varphi:\mathcal{F}(G)\rightarrow \mathcal{F}(H)$,
where $\varphi(S)=\varphi(g_1)\cdots \varphi(g_\ell)$.
If $\varphi$ is a homomorphism,
then $\varphi(S)$ is a zero-sum sequence over $H$
if and only if $\sigma(S)\in \ker(\varphi)$.
Moreover, $\varphi( \sigma(S) )=\sigma( \varphi(S) )$.

We collect some known results which will be used later.

\begin{lemma}[{\cite{[GT]}}]\label{lem:egz}
Let $n,k$ be positive integers. Every sequence of length $(k+1)n-1$ over $C_n$ has a zero-sum subsequence of length $kn$.
\end{lemma}

\begin{lemma}[{\cite[Lemma~2.2]{[GHPS]}}]\label{lem:inverse}
Let $n,k,t$ be positive integers with $2\le t<n/2+1$, and let $S$ be a sequence over $C_n$ of length $(k+1)n-t$. If $S$ contains no zero-sum subsequence of length $kn$, then there exist distinct $x,y\in C_n$ such that
\[
\mult_x(S)+\mult_y(S)\ge(k+1)n-2t+2.
\]
\end{lemma}

\begin{lemma}[Kleitman--Lemke~\cite{[KL]}]\label{lem:tiny}
Let $n$ be a positive integer. Every sequence of length $n$ over $C_n$ has a nonempty zero-sum subsequence $T$ with $\knum(T)\le1$.
\end{lemma}

For a positive integer $n$, let $\tau(n)$ denote its number of positive divisors, $\omega(n)$ its number of distinct prime divisors, and $\mult_p(n)$ the exponent of a prime $p$ in $n$. All sums indexed by $p\mid n$ are over the distinct prime divisors of $n$.

\begin{lemma}\label{lem:divisors}
Let $n$ be a positive integer and let $k=\omega(n)$.
\begin{enumerate}
\item[(i)] If $n\ge3$ is odd, then $k(\tau(n)-3)\le(n+3)/4$.
\item[(ii)] If $4\mid n$ and $k\ge4$, then $n\ge4k\tau(n)$.
\item[(iii)] If $n$ is odd and $k\ge4$, then $n\ge8k\tau(n)$.
\end{enumerate}
\end{lemma}
\begin{proof}
For every prime $p$, the ratio $p^a/(a+1)$ increases with $a\ge1$, since the ratio of consecutive terms is $p(a+1)/(a+2)>1$.

For (i), let $n=p_1^{a_1}\cdots p_k^{a_k}$ with $p_1<\cdots<p_k$ odd primes. If $k=1$, then $4(a_1-2)\le3^{a_1}+3\le n+3$. If $k=2$, it suffices to prove $8\tau(n)\le n+27$. We prove it by induction. This holds for $a_1=a_2=1$, since
$\tau(n)=4$ and $n\ge15$. Suppose that $8\tau(m)\le m+27$,
where $m$ has exactly two distinct odd
prime divisors. For any prime $p\mid m$,
we have $\tau(pm)\le3\tau(m)/2$.
Thus,
\[
8\tau(pm)\le12\tau(m)\le\frac32(m+27)\le3m+27\le pm+27.
\]
If $k\ge3$,
\[
\frac{n}{\tau(n)}\ge\frac32\cdot\frac52\cdot\frac72\cdot\left(\frac32\right)^{k-3}\ge\frac{105}{8}\left(\frac32\right)^{k-3}\ge4k,
\]
which proves (i).

For (ii), the exponent of $2$ is at least $2$, so
\[
\frac{n}{\tau(n)}\ge\frac43\cdot\frac32\cdot\frac52\cdot\frac72\cdot
\left(\frac32\right)^{k-4}
=\frac{35}{2}\left(\frac32\right)^{k-4}\ge4k.
\]
For (iii), similarly,
\[
\frac{n}{\tau(n)}\ge\frac32\cdot\frac52\cdot\frac72\cdot\frac{11}{2}\cdot
\left(\frac32\right)^{k-4}
=\frac{1155}{16}\left(\frac32\right)^{k-4}\ge8k.
\]
\end{proof}

For each $d\mid n$, let $S_d\in\mathcal F(dC_n)$
be a prescribed sequence, and let
$S=\prod_{d\mid n}S_d$.
The subsequences $S_d$ are distinguished by their
indices $d$, even if their supports intersect. For example, $S_1$ is a prescribed subsequence, not the subsequence of all terms belonging to $C_n$.

\begin{proposition}\label{prop:sequence}
Let $n$ be a positive integer, and let $S_d\in\F(dC_n)$ for every $d\mid n$. Suppose that one of the following holds:
\begin{enumerate}
\item[(i)] $n$ is odd and $\sum_{d\mid n}|S_d|=2n-1$;
\item[(ii)] $n$ is even, $\sum_{d\mid n}|S_d|=3n/2-1$, and
$\supp(S_d)\cap2dC_n=\varnothing$ if $2d\mid n$.
\end{enumerate}
Then there exist $T_d\mid S_d$ such that
\begin{equation}\label{eq:target}
\sum_{d\mid n}\sigma(T_d)=0,
\quad \sum_{d\mid n}d|T_d|=n.
\end{equation}
\end{proposition}

We now give a short proof of Theorem~\ref{thm:main}
using Proposition~\ref{prop:sequence}.

\noindent {\bf Proof of Theorem~\ref{thm:main}.}
Let $S$ be a sequence of length $2n-1$ for odd $n$, or $3n/2-1$ for even $n$ over $C_n$. For each $d\mid n$, let $S_d$ consist of its terms of order exactly $n/d$. Such terms generate $dC_n$, so $S_d\in\F(dC_n)$ and, if $2d\mid n$, no term belongs to the proper subgroup $2dC_n$. Proposition~\ref{prop:sequence} gives $T_d\mid S_d$ satisfying~\eqref{eq:target}. Then $T=\prod_{d\mid n}T_d$ is nonempty and zero-sum, and
\[
\knum(T)=\sum_{d\mid n}\frac{|T_d|}{n/d}
=\frac1n\sum_{d\mid n}d|T_d|=1.
\]
Together with the lower bounds, this proves the theorem.
\hfill$\Box$

In Proposition~\ref{prop:sequence}, the terms of $S_d$ are not
required to have order exactly $n/d$.
This is essential for the inductive argument.
If $pd\mid n$ and $R\mid S_d$ satisfies
$|R|=p$ and $\sigma(R)\in pdC_n$, then replacing
$R$ by the single term $\sigma(R)$ in $S_{pd}$
preserves the total sum and the value of
$\sum_{d\mid n}d|S_d|$, since $d|R|=pd$.
This does not require
$\ord(\sigma(R))=n/(pd)$.

\section{Odd order}

\noindent {\bf Proof of Proposition~\ref{prop:sequence}(i).}
We prove the assertion by induction on the odd integer $n$. The assertion is immediate for $n=1$. Let $n>1$, and suppose it holds for every positive odd integer $d<n$. Suppose, to the contrary, that a sequence $S=\prod_{d\mid n}S_d$ of length $2n-1$ has no subsequences satisfying~\eqref{eq:target}. Since $nC_n=\{0\}$, we have $S_n=1_{\F(C_n)}$.
For every prime $p\mid n$, let $X_p=\sum_{d\mid n,\,p\mid d}|S_d|$, and let $M=\sum_{d\mid n,\,d>1}|S_d|$.

We first prove that $X_p\le\tau(n)-3$ for every prime $p\mid n$.

Choose a prime $p\mid n$. For each $d\mid n$ with $p\nmid d$, let $\pi_p:dC_n\rightarrow dC_n/pdC_n\cong C_p$
be the canonical homomorphism. By Lemma~\ref{lem:egz}, choose a maximal family of pairwise disjoint subsequences $T_{d,1},\ldots,T_{d,a_d}\mid S_d$ with $|T_{d,i}|=p$ and $\sigma(\pi_p(T_{d,i}))=0$, and hence $\sigma(T_{d,i})\in pdC_n$ for every $i\in[1,a_d]$. Let $R_d=S_d\cdot(T_{d,1}\cdot\ldots\cdot T_{d,a_d})^{-1}$, and hence $|R_d|\le2p-2$. 

Let $\varphi:C_{n/p}\to pC_n$ be an isomorphism. For each $e\mid n/p$, set
\begin{equation}\label{eq:reduced-factors}
U_e=
\begin{cases}
\displaystyle\varphi^{-1}(S_{pe})\cdot
\prod_{i=1}^{a_e}\varphi^{-1}(\sigma(T_{e,i})),&p\nmid e,\\[5pt]
\varphi^{-1}(S_{pe}),&p\mid e.
\end{cases}
\end{equation}
Since $\varphi(eC_{n/p})=peC_n$, we have $U_e\in\F(eC_{n/p})$ and $\sum_{e\mid n/p}|U_e|=X_p+\sum_{d\mid n,\,p\nmid d}a_d$.

Suppose that $\sum_{e\mid n/p}|U_e|\ge 2(n/p)-1$.
By the induction hypothesis, there exist $V_e\mid U_e$ such that
\[
\sum_{e\mid n/p}\sigma(V_e)=0,
\quad
\sum_{e\mid n/p}e|V_e|=\frac np.
\]
Write $\varphi(V_e)=A_e\cdot\prod_{i\in I_e}\sigma(T_{e,i}),$
where $A_e\mid S_{pe}$, $I_e\subseteq[1,a_e]$ if
$p\nmid e$, and $I_e=\varnothing$ if $p\mid e$. Then $|V_e|=|A_e|+|I_e|$.
For each $d\mid n$, let
\[
T_d=
\begin{cases}
\displaystyle\prod_{i\in I_d}T_{d,i},&p\nmid d,\\[2pt]
A_{d/p},&p\mid d.
\end{cases}
\]
Since the subsequences $T_{d,1},\ldots,T_{d,a_d}$ are pairwise
disjoint, we have $T_d\mid S_d$ for every $d\mid n$. Moreover,
\[
\begin{aligned}
\sum_{d\mid n}\sigma(T_d)
&=
\sum_{e\mid n/p}
\left(\sigma(A_e)+\sum_{i\in I_e}\sigma(T_{e,i})\right)\\
&=
\varphi\left(\sum_{e\mid n/p}\sigma(V_e)\right)=0,
\end{aligned}
\]
and 
\[
\begin{aligned}
\sum_{d\mid n}d|T_d|
&=
\sum_{\substack{d\mid n\\p\mid d}}d|A_{d/p}|
+
\sum_{\substack{d\mid n\\p\nmid d}}pd|I_d|\\
&=
p\sum_{e\mid n/p}e\bigl(|A_e|+|I_e|\bigr)
=
p\sum_{e\mid n/p}e|V_e|
=n.
\end{aligned}
\]
This contradicts our assumption. Hence,
\begin{equation}\label{eq:compression}
X_p+\sum_{\substack{d\mid n\\p\nmid d}}a_d\le\frac{2n}{p}-2.
\end{equation}

Since $\sum_{d\mid n,\,p\nmid d}|S_d|=2n-1-X_p$, multiplying~\eqref{eq:compression} by $p$ yields
\[
(p-1)X_p+2p-1
\le\sum_{\substack{d\mid n\\p\nmid d}}|R_d|.
\]
Since precisely $\tau(n)/(\mult_p(n)+1)$ positive divisors of $n$ are not divisible by $p$, we have
\[
\sum_{\substack{d\mid n\\p\nmid d}}|R_d|
\le\frac{(2p-2)\tau(n)}{\mult_p(n)+1}\le(p-1)\tau(n).
\]
It follows that $X_p\le\tau(n)-(2p-1)/(p-1)<\tau(n)-2$, and therefore $X_p\le\tau(n)-3$.

By the definition of $X_p$ and Lemma~\ref{lem:divisors}(i), we have
\begin{equation}\label{eq:M}
M
=\sum_{\substack{d\mid n\\d>1}}|S_d|
\le
\sum_{\substack{d\mid n\\d>1}}\omega(d)|S_d|
=\sum_{p\mid n}X_p
\le\omega(n)(\tau(n)-3)
\le\frac{n+3}{4}.
\end{equation}

If $n$ is prime, the estimate $X_n\le\tau(n)-3=-1$ is already a contradiction. We may therefore assume that $n$ is composite, so $n\ge9$.

Since the sequence $S_1$ of length $2n-1-M$ has no zero-sum subsequence of length $n$, by Lemma~\ref{lem:egz}, we have $M\ge1$. Taking $t=M+1$, we have $2\le t<n/2+1$ by~\eqref{eq:M}. By Lemma~\ref{lem:inverse}, with $k=1$, there exist distinct $x,y\in C_n$ such that $\mult_x(S_1)+\mult_y(S_1)\ge2n-2M$.

Since $S_1$ contains no zero-sum subsequence of length $n$ and $ng=0$ for every $g\in C_n$, we have $\mult_g(S_1)\le n-1$ for every $g\in C_n$. Write $\mult_x(S_1)=n-1-\alpha$ and $\mult_y(S_1)=n-1-\beta$, where $\alpha,\beta\ge0$. Then
\[
\alpha+\beta\le2M-2\le\frac{n-1}{2}.
\]
We claim that $u=x-y$ generates $C_n$. Otherwise, $e=\ord(u)$ is a proper divisor of $n$, so $e\le n/3$. For every integer $a$ with $\beta+1\le a\le n-1-\alpha$, the sequence $x^{a}\cdot y^{n-a}$ divides $S_1$ and has sum $au$. This interval contains $n-1-\alpha-\beta\ge(n-1)/2\ge e$ consecutive integers, and hence contains a multiple of $e$. The corresponding subsequence is zero-sum and has length $n$, a contradiction.

Let $r=\alpha+\beta+1$. Choose $g_1,\ldots,g_r\in C_n$
and divisors $d_1,\ldots,d_r$ of $n$ such that, for every $d\mid n$,
\[
\prod_{\substack{1\le j\le r\\d_j=d}}g_j
=
\begin{cases}
S_1\bigl(x^{n-1-\alpha}\cdot y^{n-1-\beta}\bigr)^{-1},
&d=1,\\[2pt]
S_d,&d>1.
\end{cases}
\]
Then
$S=x^{n-1-\alpha}\cdot y^{n-1-\beta}
  \cdot g_1\cdot\ldots\cdot g_r,$
where $g_j\in d_jC_n$ for every $j\in[1,r]$, and
$g_j\notin\{x,y\}$ whenever $d_j=1$. Since $u=x-y$ generates $C_n$ and $d_jy-g_j\in d_jC_n$,
there is a unique
$D_j\in\{0,d_j,2d_j,\ldots,n-d_j\}$ such that
$g_j=d_jy-D_ju$.

Consider the subsequence $T=x^a\cdot y^b\cdot\prod_{j\in J}g_j\mid S$, where
$J\subseteq[1,r]$, $0\le a\le n-1-\alpha$ and $0\le b\le n-1-\beta$,
with the factorization $T_1=x^ay^b\prod_{j\in J,\,d_j=1}g_j$ and
$T_d=\prod_{j\in J,\,d_j=d}g_j$ for $d>1$. Then
\[
\sum_{d\mid n}d|T_d|=a+b+\sum_{j\in J}d_j ,
\]
and,
\[
\sigma(T)
=\left(a+b+\sum_{j\in J}d_j\right)y
+\left(a-\sum_{j\in J}D_j\right)u.
\]
Hence, $T$ satisfies~\eqref{eq:target} as soon as
$a+b+\sum_{j\in J}d_j=n$ and $a\equiv\sum_{j\in J}D_j\pmod n$.

We shall show that these two conditions hold for some
$J\subseteq[1,r]$, $0\le a\le n-1-\alpha$, and
$0\le b\le n-1-\beta$.

Suppose first that there exists $j\in[1,r]$ such that
$D_j\le n-1-\alpha$ and $D_j+d_j\ge\beta+1$.
Take $J=\{j\}$, $a=D_j$, and $b=n-d_j-D_j$. Then these two conditions hold.

Otherwise, for every $j\in[1,r]$, either
$D_j+d_j\le\beta$ or $n-D_j\le\alpha$.
These two inequalities cannot both hold, since their sum
would give $n+d_j\le\alpha+\beta\le(n-1)/2$.
Thus
\[
Y=\{j\in[1,r]:D_j+d_j\le\beta\},
\quad
X=\{j\in[1,r]:n-D_j\le\alpha\}
\]
partition $[1,r]$.

For $j\in Y$, we have $D_j+d_j\ge2$. Otherwise, $d_j=1$
and $D_j=0$, which gives $g_j=y$, contradicting
$g_j\notin\{x,y\}$ whenever $d_j=1$.
Similarly, for $j\in X$, we have $n-D_j\ge2$. Otherwise,
$D_j=n-1$ and $d_j=1$, which gives $g_j=x$, contradicting
$g_j\notin\{x,y\}$ whenever $d_j=1$. Consequently,
\[
\sum_{j\in Y}(D_j+d_j)+\sum_{j\in X}(n-D_j)
\ge2r>\alpha+\beta.
\]

Suppose that $\sum_{j\in Y}(D_j+d_j)\ge\beta+1$.
Choose a minimal subset $J\subseteq Y$ such that
$\sum_{j\in J}(D_j+d_j)\ge\beta+1$. Since
$D_j+d_j\le\beta$ for $j\in Y$, minimality gives
$\beta+1\le\sum_{j\in J}(D_j+d_j)\le2\beta.$
Take $a=\sum_{j\in J}D_j$ and
$b=n-\sum_{j\in J}(D_j+d_j)$. Since
$\alpha+\beta\le(n-1)/2$, we obtain
$0\le a\le2\beta\le n-1-\alpha$, 
$0\le n-2\beta\le b\le n-1-\beta.$
Hence, these two conditions hold.

It remains to consider
$\sum_{j\in Y}(D_j+d_j)\le\beta$.
Then $\sum_{j\in X}(n-D_j)\ge\alpha+1$.
Choose a minimal subset $J\subseteq X$ such that
$\sum_{j\in J}(n-D_j)\ge\alpha+1$. Since
$n-D_j\le\alpha$ for $j\in X$, we have
$\alpha+1\le\sum_{j\in J}(n-D_j)\le2\alpha.$
Take $a=n-\sum_{j\in J}(n-D_j)$ and
$b=\sum_{j\in J}(n-D_j-d_j)$.
Since $D_j\le n-d_j$ and $\alpha+\beta\le(n-1)/2$, we obtain
$0\le n-2\alpha\le a\le n-1-\alpha$, $0\le b\le2\alpha\le n-1-\beta$.
Thus, these two conditions also hold in this situation.

In each situation, the resulting subsequences $T_d\mid S_d$
satisfy~\eqref{eq:target}, contrary to our assumption.
This completes the induction.
\hfill$\Box$

\section{Even order}

For $T=\prod_{d\mid n}T_d$ with $T_d\mid S_d$, let $\wt(T)=\sum_{d\mid n}d|T_d|$, where the terms of each $T_d$ are taken from the prescribed sequence $S_d$. Then $|T|\le\wt(T)$, and $\wt(UV)=\wt(U)+\wt(V)$ for disjoint subsequences $U,V$ with these factorizations. Thus~\eqref{eq:target} is equivalent to $\sigma(T)=0$ and $\wt(T)=n$.

\begin{lemma}\label{lem:completion}
Let $n=2m$ be a positive even integer, and let $S=\prod_{d\mid n}S_d$, where $S_d\in\F(dC_n)$, $|S|=3m-1$, and $\supp(S_d)\cap2dC_n=\varnothing$ whenever $2d\mid n$. Suppose that $x\in C_n\setminus2C_n$ and $P\mid S$ satisfy
$m-1\le\wt(P)\le m$, and that, for every $j\in[0,\wt(P)]$, there exists $Q\mid P$ such that
$\wt(Q)=j$ and $\sigma(Q)=jx$.
Then $S$ contains a zero-sum subsequence $T$ with $\wt(T)=n$.
\end{lemma}

\begin{proof}
Choose $U\mid SP^{-1}$ satisfying $\sigma(U)=\wt(U)x$ and $\wt(U)\le n-1-\wt(P)$
such that $\wt(U)$ is maximal, where the empty sequence is admitted.
Since $|P|\le\wt(P)$ and $|U|\le\wt(U)$, we have
\[
|S(PU)^{-1}|
\ge3m-1-\wt(P)-\wt(U)
\ge3m-1-(n-1)=m.
\]
Let $R=g_1\cdot\ldots\cdot g_m$ be a subsequence of $S(PU)^{-1}$, where, for each $i\in[1,m]$, $g_i$ is a term of $S_{d_i}$ with $d_i\mid n$ .

If $d_i$ is even, then $g_i,d_ix\in d_iC_n\subseteq2C_n$.
If $d_i$ is odd, then $2d_i\mid n$ and
$d_iC_n\cap2C_n=2d_iC_n$.
Since $\supp(S_{d_i})\cap2d_iC_n=\varnothing$, we have $g_i\in d_iC_n\setminus2d_iC_n$ and hence $g_i\notin2C_n$.
Also, $d_ix\notin2C_n$, since $x\notin2C_n$ and $d_i$ is odd.
Thus, we have $g_i-d_ix\in d_iC_n\cap2C_n$ for any $i\in[1,m]$, and hence $\ord(g_i-d_ix)\mid n/d_i$.

Consider the sequence
\[
R'=\prod_{i=1}^{m}(g_i-d_ix)\in\F(2C_n).
\]
Since $|R'|=m=|2C_n|$, by Lemma~\ref{lem:tiny}, there exists a nonempty subsequence $V'=\prod_{i\in I}(g_i-d_ix)\mid R'$ such that $\sigma(V')=0$ and $\knum(V')=\sum_{i\in I}\frac1{\ord(g_i-d_ix)}\le1$, where $\varnothing\ne I\subseteq[1,m]$.
Let $V=\prod_{i\in I}g_i\mid R$. Then $\sigma(V)=\sum_{i\in I}g_i=\sum_{i\in I}d_ix=\wt(V)x$ and $1\le\wt(V)=\sum_{i\in I}d_i\le n\sum_{i\in I}\frac1{\ord(g_i-d_ix)}\le n$.
Thus, $U$ and $V$ are disjoint, and
$\sigma(UV)=\wt(UV)x$.
The maximality of $\wt(U)$ therefore implies $\wt(U)+\wt(V)\ge n-\wt(P).$
If $\wt(V)\le n-\wt(P)-1$, we take $V^*=UV$. Since $\wt(P)\ge m-1$, we have
\[
n-\wt(P)\le\wt(V^*)
=\wt(U)+\wt(V)
\le2\bigl(n-1-\wt(P)\bigr)\le n.
\]
If $\wt(V)\ge n-\wt(P)$, we take $V^*=V$.
Thus, we have $V^*\mid SP^{-1}$, $n-\wt(P)\le\wt(V^*)\le n$ and $\sigma(V^*)=\wt(V^*)x$.
Since $0\le n-\wt(V^*)\le\wt(P)$, the hypothesis gives
$Q\mid P$ such that $\wt(Q)=n-\wt(V^*)$ and $\sigma(Q)=\bigl(n-\wt(V^*)\bigr)x$. Hence, $T=QV^*$ is the required subsequence.
\end{proof}

\noindent {\bf Proof of Proposition~\ref{prop:sequence}(ii).}
We prove the assertion by induction on the even integer $n$. The assertion is immediate for $n=2$. Let $n=2m>2$, and suppose it holds for every positive even integer $d<n$. Suppose, to the contrary, that a sequence $S=\prod_{d\mid n}S_d$ of length $3n/2-1$ has no subsequences satisfying~\eqref{eq:target}. Since $nC_n=\{0\}$, we have $S_n=1_{\F(C_n)}$. For every prime $p\mid n$, let $X_p=\sum_{d\mid n,\,p\mid d}|S_d|$, and let $M=\sum_{d\mid n,\,d>1}|S_d|$.

We first prove that
\begin{equation}\label{eq:Xp-cap}
X_p\le\frac{2\tau(n)}{\mult_p(n)+1}-3,
\end{equation}
for every prime $p\mid n$ such that $n/p$ is even. Choose a prime $p\mid n$ such that $n/p$ is even, and let
$d\mid n$ with $p\nmid d$.

Suppose that $p$ is odd.
If $|S_d|\ge2p-1$, let $R\mid S_d$ with $|R|=2p-1$.
Applying Lemma~\ref{lem:egz} to the image of $R$ in
$dC_n/pdC_n\cong C_p$, we obtain a subsequence $T\mid R$ of length $p$ satisfying $\sigma(T)\in pdC_n$.
When $2pd\mid n$, we have
$\supp(S_d)\subseteq dC_n\setminus2dC_n$.
Let $\pi_2:dC_n\rightarrow dC_n/2dC_n\cong C_2$
be the canonical homomorphism.
Thus $\pi_2(T)=h^p$, where $h$ is the unique nonzero element of $dC_n/2dC_n$.
Since $p$ is odd and $\ord(h)=2$, we have $\pi_2(\sigma(T))=\sigma(\pi_2(T))=ph=h\ne0.$
Thus, $\sigma(T)\in pdC_n\setminus2pdC_n$.

Suppose now that $p=2$.
Then $4\mid n$ and $d$ is odd, so $4d\mid n$.
Let $\pi_4:dC_n\rightarrow dC_n/4dC_n\cong C_4$
be the canonical homomorphism. If $|S_d|\ge3$, 
let $R=g_1g_2g_3\mid S_d$, where $g_i\in dC_n\setminus2dC_n$ for every $i\in[1,3]$, so $\ord(\pi_4(g_i))=4$.
Since $dC_n/4dC_n\cong C_4$ has exactly two elements
of order $4$, there exist distinct $i,j\in[1,3]$
such that $\pi_4(g_i)=\pi_4(g_j)$.
For $T=g_ig_j\mid R$, we have $|T|=2$ and $\pi_4(\sigma(T))
=\pi_4(g_i)+\pi_4(g_j)
=2\pi_4(g_i)$.
It follows that $\sigma(T)\in2dC_n\setminus4dC_n.$

Consequently, as in the odd-order proof, there is a factorization $S_d=T_{d,1}\cdot\ldots\cdot T_{d,a_d}\cdot R_d$ such that $|T_{d,i}|=p$ for $i\in[1,a_d]$, $|R_d|\le2p-2$, and
\[
\sigma(T_{d,i})\in pdC_n,\quad
\sigma(T_{d,i})\notin2pdC_n\ \text{if }2pd\mid n.
\]
Let $\varphi:C_{n/p}\to pC_n$ be an isomorphism and define $U_e$ by~\eqref{eq:reduced-factors}. Since $\varphi(2eC_{n/p})=2peC_n$, we have $U_e\in\mathcal F(eC_{n/p})$ for every $e\mid n/p$, and $\supp(U_e)\cap2eC_{n/p}=\varnothing$ whenever $2e\mid n/p$.
If $\sum_{e\mid n/p}|U_e|\ge 3n/(2p)-1$, the induction hypothesis and the same construction of $T_d$ as in the odd-order proof would give~\eqref{eq:target}. Therefore
\[
X_p+\sum_{\substack{d\mid n\\p\nmid d}}a_d
=\sum_{e\mid n/p}|U_e|\le\frac{3m}{p}-2.
\]
Since $|S_d|=pa_d+|R_d|$ and $\sum_{d\mid n,\,p\nmid d}|S_d|=3m-1-X_p$, we obtain
\begin{equation}\label{eq:even-count}
(p-1)X_p+2p-1\le\sum_{\substack{d\mid n\\p\nmid d}}|R_d|
\le\frac{(2p-2)\tau(n)}{\mult_p(n)+1}.
\end{equation}
Thus, $X_p\le\frac{2\tau(n)}{\mult_p(n)+1}-3.$

We claim that
$\h(S_1)\le m-2$. For every $x\in C_n$, if $x^{m-1}\mid S_1$, then $x\notin2C_n$. Taking $P=x^{m-1}$ in Lemma~\ref{lem:completion} gives a contradiction, since $Q=x^j$ satisfies $\wt(Q)=j$ and $\sigma(Q)=jx$ for every $j\in[0,m-1]$.

Fix $g_0\in C_n\setminus2C_n$.
Let $S_1=g_1\cdot\ldots\cdot g_{|S_1|}$ and let
\[
W=\prod_{i=1}^{|S_1|}(g_i-g_0)\in\F(2C_n)\cong \F(C_m).
\]
If $W$ contains a zero-sum subsequence
$W'=\prod_{i\in I}(g_i-g_0)$, where $I\subseteq[1,|S_1|]$ and $|I|=2m$, then the sequence
$T=\prod_{i\in I}g_i\mid S_1$ satisfies \eqref{eq:target}, a contradiction.
Thus, $W$ has no zero-sum subsequence of length $2m$,
and Lemma~\ref{lem:egz} gives $|W|\le3m-2$.
If $|W|\ge\lfloor5m/2\rfloor$, then $t=3m-|W|$ satisfies
$2\le t\le\lceil m/2\rceil<m/2+1$.
By Lemma~\ref{lem:inverse}, with $k=2$,
there exist distinct $h_1,h_2\in2C_n\cong C_m$ such that
\[
\mult_{h_1}(W)+\mult_{h_2}(W)
\ge3m-2t+2\ge2m+1.
\]
Since $\mult_{h_i}(W)=\mult_{h_i+g_0}(S_1)$ for $i\in\{1,2\}$,
this contradicts $\h(S_1)\le m-2$.
Therefore $|S_1|\le\lfloor5m/2\rfloor-1$,
which proves $M\ge\left\lceil\frac m2\right\rceil.$

Suppose that $m$ is even. Then $n/p$ is even for every prime $p\mid n$. By the definition of $M$,
\[
M\le\sum_{\substack{d\mid n\\d>1}}\omega(d)|S_d|=\sum_{p\mid n}X_p.
\]
Moreover, for each prime $p\mid n$,
\[
\sum_{\substack{d\mid n\\p\nmid d}}|R_d|
\le |R_1|+\sum_{\substack{d\mid n\\d>1,\ p\nmid d}}|S_d|
\le2p-2+M-X_p.
\]
Together with~\eqref{eq:even-count}, this gives
$pX_p\le M-1$ for every prime $p\mid n.$
If $\omega(n)\le2$, then $M\le(5/6)(M-1)$, a contradiction. If $\omega(n)=3$, set $n=2^{\alpha}p^{\beta}q^{\lambda}$, where $\alpha\ge2$, $\beta,\lambda\ge1$, and $3\le p<q$. By~\eqref{eq:Xp-cap}, we have
\[
M\le X_2+X_p+X_q\le2(\beta+1)(\lambda+1)-3+\frac8{15}(M-1),
\]
and hence
\[
7M\le30(\beta+1)(\lambda+1)-53.
\]
Since $\beta+\lambda\le2\beta \lambda$, $p^\beta\ge3\beta$, and $q^\lambda\ge5\lambda$, we obtain
\[
7M\le90\beta\lambda-23<105\beta \lambda\le7p^\beta q^\lambda\le\frac{7n}{4},
\]
contrary to $M\ge\left\lceil\frac m2\right\rceil\ge n/4$. If $\omega(n)\ge4$, Lemma~\ref{lem:divisors}(ii) and~\eqref{eq:Xp-cap} yield
\[
M\le\sum_{p\mid n}\left(\frac{2\tau(n)}{\mult_p(n)+1}-3\right)
\le \omega(n)(\tau(n)-3)<n/4,
\]
again a contradiction.

It remains to consider $n=2m$ with $m\ge3$ odd. Then $n/p$ is even for every prime $p\mid m$. Let $b=|S_2|$. Then
\[
M-b=\sum_{\substack{d\mid m\\d>1}}(|S_d|+|S_{2d}|)
\le\sum_{\substack{d\mid m\\d>1}}\omega(d)(|S_d|+|S_{2d}|)
=\sum_{p\mid m}X_p.
\]
For each prime $p\mid m$, we also have
\[
\sum_{\substack{d\mid n\\p\nmid d}}|R_d|
\le |R_1|+|R_2|+
\sum_{\substack{d\mid n\\d\notin\{1,2\},\ p\nmid d}}|S_d|
\le4p-4+M-b-X_p.
\]
Together with~\eqref{eq:even-count}, this gives
\begin{equation}\label{eq:Xp-linear-two}
pX_p\le M-b+2p-3,
\end{equation}
for every prime $p\mid m$.
We claim that
\begin{equation}\label{eq:even-M}
M-b\le\frac{m+1}{4}.
\end{equation}
For $\omega(m)=1$,~\eqref{eq:Xp-cap} gives $M-b\le1$. For $\omega(m)=2$ and $M-b\ge3$,~\eqref{eq:Xp-linear-two} gives
\[
M-b\le\frac8{15}(M-b-3)+4,
\]
so $M-b\le5$. If $M-b=5$, then $X_p\le\lfloor(2p+2)/p\rfloor=2$ for both primes, giving $M-b\le4$, a contradiction. Thus $M-b\le4\le(m+1)/4$. For $\omega(m)=3$ and $M-b\ge3$, by \eqref{eq:Xp-linear-two}, we have
\[
M-b\le\left(\frac13+\frac15+\frac17\right)(M-b-3)+6
=\frac{71}{105}(M-b-3)+6,
\]
so $M-b\le12<(m+1)/4$, since $m\ge105$. 
If $\omega(m)\in\{2,3\}$ and $M-b\le2$, then~\eqref{eq:even-M} holds
trivially. For $\omega(m)\ge4$, Lemma~\ref{lem:divisors}(iii) and~\eqref{eq:Xp-cap} give
\[
M-b\le\sum_{p\mid m}\left(\frac{4\tau(m)}{\mult_p(m)+1}-3\right)
\le\omega(m)(2\tau(m)-3)<m/4.
\]
This proves the claim.

Since $m$ is odd, $C_n\cong C_2\oplus C_m$. Let $\pi:C_n\to C_m$ be the projection. Then $\supp(S_1)\subseteq\{1\}\times C_m$ and $\supp(S_2)\subseteq\{0\}\times C_m$. By \eqref{eq:even-M} and $M\ge\left\lceil\frac m2\right\rceil$, we have
\begin{equation}\label{eq:even-bM}
b=M-(M-b)\ge M-\frac{m+1}{4}\ge\frac M2.
\end{equation}

Let $e$ be a generator of $C_m$. Pairing equal terms of $S_1$ as far as possible, we may write
\[
S_1=x_1^2\cdot\ldots\cdot x_r^2\cdot
(1,z_1e)\cdot\ldots\cdot(1,z_te),
\]
where $0\le z_1<\cdots<z_t\le m-1$ and $t\le m$. Set $P_i=x_i^2$ for $i\in[1,r]$ and
$P_{r+j}=(1,z_{2j-1}e)\cdot(1,z_{2j}e)$ for
$j\in[1,\lfloor t/2\rfloor]$. If $t$ is odd, let
$R'=(1,z_te)$. Otherwise, let $R'$ be the empty sequence. Then
$S_1=P_1\cdot\ldots\cdot P_s\cdot R'$, where $|P_i|=2$, $|R'|\le1$ and $s=r+\left\lfloor\frac t2\right\rfloor
=\left\lfloor\frac{|S_1|}{2}\right\rfloor$.
Let
\[
V=\pi(S_2)\cdot\prod_{i=1}^s\pi(\sigma(P_i))\in\F(C_m),
\]
and let $\ell=|V|=b+\left\lfloor\frac{|S_1|}{2}\right\rfloor.$

If $V$ has a zero-sum subsequence of length $m$, there exist $A\mid S_2$ and
$I\subseteq[1,s]$ such that $|A|+|I|=m$ and
$\pi(\sigma(A)+\sum_{i\in I}\sigma(P_i))=0$.
Since $\sigma(A)+\sum_{i\in I}\sigma(P_i)\in2C_n$ and
$\pi$ is injective on $2C_n$, it follows that
$T_1=\prod_{i\in I}P_i$ and $T_2=A$, with all other
$T_d$ empty, satisfy~\eqref{eq:target}, a contradiction.
By Lemma~\ref{lem:egz}, we have $\ell\le2m-2$.

If $m=3$, we have $\h(S_1)\le m-2=1$, and hence $|S_1|\le3$. By~\eqref{eq:even-M}, we have $M-b\le1$.
Since $M=8-|S_1|$, $\ell\ge7-\lceil|S_1|/2\rceil\ge5>2m-2$, a contradiction.
Thus $m\ge5$, and~\eqref{eq:even-M} implies
$M-b\le\lfloor(m+1)/4\rfloor\le(m-3)/2$.

For every $h\in\supp(\pi(S_1))$, we claim that
\begin{equation}\label{eq:pair-completion}
\mult_{2h}(V)\le\frac{m-3}{2}.
\end{equation}
Assume that there exists $h\in\supp(\pi(S_1))$ satisfying
$\mult_{2h}(V)\ge(m-1)/2$, and let $x=(1,h)$.
By the definition of $V$, there exist pairwise disjoint
subsequences $W_1,\ldots,W_{(m-1)/2}$ of $S$ such that
$\wt(W_i)=2$ and $\sigma(W_i)=2x$ for every $i\in[1,(m-1)/2]$.
Each $W_i$ is either a term of $S_2$ or one of the pairs
$P_j$ for some $j\in[1,s]$.

If $x\nmid W_1\cdots W_{(m-1)/2}$, let $P=x\cdot W_1\cdots W_{(m-1)/2}$. Then
$\wt(P)=m$, and for $j\in[0,m]$ we take $Q=x^{\varepsilon}\prod_{i\in I}W_i$
with $\varepsilon\in\{0,1\}$, $\varepsilon\equiv j\pmod 2$ and
$|I|=(j-\varepsilon)/2$.
If $x\mid W_1\cdots W_{(m-1)/2}$, there exists $i_0\in[1,(m-1)/2]$ such that $x\mid W_{i_0}$. Since $\sigma(W_{i_0})=2x$, we have $W_{i_0}=x^2$. Let $P=W_1\cdots W_{(m-1)/2}$, so $\wt(P)=m-1$.
For even $j$, we take $Q=\prod_{i\in I}W_i$ with $|I|=j/2$. For odd $j$, we take
$Q=x\cdot\prod_{i\in I}W_i$ with $I\subseteq[1,\frac{m-1}{2}]\setminus\{i_0\}$
and $|I|=\frac{j-1}{2}\le\frac{m-3}{2}$.
In either case, we have $\wt(Q)=j$ and $\sigma(Q)=jx$. This contradicts
Lemma~\ref{lem:completion} and proves~\eqref{eq:pair-completion}.

By~\eqref{eq:even-bM} and $|S_1|=3m-1-M$, we have
\[
2\ell\ge 2b+|S_1|-1
=3m-2+2b-M\ge 3m-2.
\]
Since $m$ is odd, $\ell\ge (3m-1)/2$.
Hence,
\[
2\le 2m-\ell\le (m+1)/2<m/2+1.
\]
Applying Lemma~\ref{lem:inverse} to $V$ with $k=1$, there exist distinct $u,v\in C_m$ satisfying
\begin{equation}\label{eq:two-pair-values}
\mult_u(V)+\mult_v(V)
\ge 2\ell-2m+2
\ge m+2b-M.
\end{equation}

For every $w\in 2\supp(\pi(S_1))=\{2h:h\in\supp(\pi(S_1))\}$, by~\eqref{eq:pair-completion}, we have $\mult_w(V)\le\frac{m-3}{2}$.
Since $1\le z_1+z_2<z_3+z_4<\cdots\le2m-3$, and $\ord(e)=m$, each element of $C_m$ occurs at most twice among the sums $(z_{2j-1}+z_{2j})e$.
Thus, for every $w\notin 2\supp(\pi(S_1))$, we have
\[
\mult_w(V)
\le \mult_w(\pi(S_2))+2.
\]

If $u,v\notin 2\supp(\pi(S_1))$ are distinct, then we have
\[
\mult_u(V)+\mult_v(V)
\le \mult_u(\pi(S_2))+\mult_v(\pi(S_2))+\lfloor t/2\rfloor\le b+\frac{m-1}{2}.
\]
If $u,v\in 2\supp(\pi(S_1))$ are distinct, then $\mult_u(V)+\mult_v(V)\le m-3.$
If exactly one of $u,v$ belongs to
$2\supp(\pi(S_1))$, then
\[
\mult_u(V)+\mult_v(V)
\le \frac{m-3}{2}+b+2
=b+\frac{m+1}{2}.
\]
Consequently, in all cases,
\[
\mult_u(V)+\mult_v(V)
\le \max\left\{
m-3,\ b+\frac{m+1}{2}
\right\}
<m+2b-M,
\]
where the last inequality follows from \eqref{eq:even-M} and \eqref{eq:even-bM}.
This contradicts~\eqref{eq:two-pair-values}
and completes the induction.
\hfill$\Box$

%\section*{Acknowledgments}


\begin{thebibliography}{9}

\setlength{\itemsep}{2pt}

\bibitem{[EGZ]} P.~Erd\H{o}s, A.~Ginzburg and A.~Ziv, \emph{Theorem in the additive number theory}, Bull. Res. Council Israel \textbf{10F} (1961), 41--43.

\bibitem{[GHPS]} W.~Gao, D.~Han, J.~Peng and F.~Sun, \emph{On zero-sum subsequences of length $k\exp(G)$}, J. Combin. Theory Ser. A \textbf{125} (2014), 240--253.

\bibitem{[GHJLW]} W.~Gao, W.~Hui, X.~Jiang, Y.~Li and X.~Wang, \emph{On zero-sum subsequences of cross number $1$}, Acta Arith. \textbf{219} (2025), 347--363.

\bibitem{[GT]}
W.~Gao and R.~Thangadurai,
\emph{On zero-sum sequences of prescribed length},
Aequationes Math. \textbf{72} (2006), 201--212.

\bibitem{[GH2006]}
A. Geroldinger and F. Halter-Koch,
\emph{Non-unique factorizations. Algebraic, Combinatorial and Analytic Theory},
Pure Appl. Math., vol. 278, Chapman $\&$ Hall/CRC, 2006.

\bibitem{[KL]} P.~Lemke and D.~Kleitman, \emph{An addition theorem on the integers modulo $n$}, J. Number Theory \textbf{31} (1989), 335--345.

\end{thebibliography}
\end{document}